\documentclass[letterpaper, 10 pt, conference]{ieeeconf}

\usepackage{amssymb,amsmath,color}
\usepackage{amsfonts,mathtools}
\usepackage{graphicx,bm}
\usepackage{caption}
\mathtoolsset{showonlyrefs}
\usepackage{dsfont}
\IEEEoverridecommandlockouts                              
\usepackage{cite}                                                          
\newcommand*{\QEDB}{\hfill\ensuremath{\square}}%

\newtheorem{assumption}{Assumption}
\newtheorem{definition}{Definition}

\newtheorem{lemma}{Lemma}
\newtheorem{theorem}{Theorem}

\newtheorem{remark}{Remark}

\begin{document}

\title{\bf Averaging in a Class of Stochastic Hybrid Dynamical Systems with Time-Varying Flow Maps}

\author{Jorge I. Poveda \thanks{Jorge I. Poveda is with the ECE Department of the Jacobs School of Engineering, University of California, San Diego, La Jolla, CA, USA. Email:\{poveda@ucsd.edu.\}.}\thanks{Research supported in part by AFOSR grant FA9550-22-1-0211, and by NSF grant ECCS CAREER 2305756.}}

\maketitle
\begin{abstract}
We present stability and recurrence results for a class of stochastic hybrid dynamical systems with oscillating flow maps. These results are developed by introducing averaging tools that parallel those already existing for ordinary differential equations and deterministic hybrid dynamical systems. Such tools can be used to examine the stability properties of the original dynamics based on the properties of a simpler dynamical system constructed from the average of the original oscillating vector field. In this work, we focus on a class of systems for which \emph{global} stability and recurrence results are achievable under suitable smoothness assumptions on the dynamics. By studying the average stochastic hybrid dynamics using Lyapunov-Foster functions, we derive similar stability and recurrence results for the original stochastic hybrid system.
\end{abstract}
\section{INTRODUCTION}
During the last century, averaging theory \cite{Krylov1, Bogoliubov1, HaleODE80} has played a significant role in synthesizing and analyzing nonlinear controllers \cite{KhalilBook}, dynamics estimation \cite{Sastry:89}, model-free optimization algorithms \cite{KrsticBookESC, PoTe17Auto}, and even in the analysis of quantum systems \cite[Sec. 1.2]{AveragingQuantumSystems}. For ordinary differential equations (ODEs) and deterministic hybrid dynamical systems (HDS) \cite{bookHDS}, averaging tools were extensively studied in \cite{TeelMoreauNesic2003} and in \cite{Wang:12_Automatica, PovedaNaliAuto20, AveragingHybridISS}.

Loosely speaking, such results rely on establishing two key technical properties:

(a) The solutions of the original time-varying system and those of its average system remain ``close" to each other on compact time domains. This property is typically obtained through a suitable change of variables that transforms the original dynamics into a perturbed version of the average dynamics, and as shown in \cite{TeelNesicAveraging, Wang:12_Automatica}, this approach can also be used in HDS using constructions similar to those developed for ODEs in \cite{Bogoliubov1, HaleODE80}.

(b) The stability properties of the average dynamics are robust concerning small perturbations. For deterministic HDS, these robustness properties can be directly obtained under suitable regularity conditions on the data of the dynamics (see \cite[Ch.7]{bookHDS}). 

The combination of both properties leads to uniform convergence results via $\mathcal{K}\mathcal{L}$ bounds, from arbitrarily large (but fixed) compact sets to arbitrarily small (but fixed) neighborhoods of the sets of interest. In the literature, this property is known as ``semi-global practical asymptotic stability" (SGPAS) \cite{TeelMoreauNesic2003}.

During the last decade, the development of averaging results for deterministic HDS has opened the door to new algorithms and controllers that are capable of overcoming some of the fundamental limitations of conventional techniques. Examples include PWM implementations of logic-based controllers \cite{TeelNesicAveraging}, hybrid extremum-seeking control \cite{PoTe17Auto, Kutadinata_Moase}, and hybrid vibrational control \cite{kapitzaochoa}.

However, when stochastic uncertainty (as opposed to worst-case uncertainty) is introduced into HDS, averaging techniques for stability analysis are largely lacking, and existing results predominantly focus on impulsive systems \cite{tsai1997averaging} or stochastic differential equations \cite{averagingODEdiffusions}. Nonetheless, as thoroughly discussed in \cite{hespanha2005model, Teel:14_AutomaticaSurvey}, stochastic hybrid dynamical systems can be employed to capture various phenomena that arise in practical engineering applications. These applications encompass systems with stochastic communication networks \cite{hespanha2005model}, automation processes in uncertain environments \cite{Cassandras_Lygeros_2010}, and energy systems with renewable generation \cite[Sec. III.2.d]{ochoa2023control}. Such applications would greatly benefit from novel algorithms and controllers for stochastic systems capable of addressing some of the limitations of existing techniques, while still relying on averaging theory to certify stability and robustness properties.

In this paper, we take a step forward in the aforementioned direction by introducing averaging tools for a class of stochastic hybrid dynamical systems (SHDS) that incorporate nonlinear flow maps and jump maps \cite{Teel_ANRC, Teel:14_AutomaticaSurvey}. In our case, the SHDS generates a well-defined average SHDS induced by a fast-varying state in the flow map of the original dynamics. The jump map introduces stochasticity through random inputs modeled by a sequence of independent and identically distributed (i.i.d.) random variables. For simplicity, we study the stability properties of the main state of the system with respect to the origin. However, the system can include additional auxiliary states, common in hybrid systems for modeling timers, logic states, counters, etc. For these states, stability properties are typically studied with respect to compact sets.

Our main result, presented in Theorem 1, demonstrates that provided the average SHDS satisfies appropriate regularity and stochastic stability properties, as certified by a quadratic-like Lyapunov-Foster function, uniform global recurrence can be established for the original SHDS whenever the time variation of the oscillating signals is sufficiently high. Furthermore, the recurrent set can be made arbitrarily close to the origin by increasing the time scale separation in the system. In cases where the flow map of the main state does not explicitly depend on the auxiliary states, or when such auxiliary states remain constant during flows (e.g., as in switched systems), we obtain uniform global exponential stability in expectation for the original SHDS. Our results are established by leveraging Lyapunov-Foster conditions introduced in the seminal paper \cite{Teel_ANRC}, as well as classic constructions studied in \cite{Bogoliubov1, HaleODE80, Sastry:89} for ODEs and in \cite{TeelNesicAveraging} for deterministic hybrid systems.

The rest of this paper is organized as follows: Section \ref{sec:preli} introduces the notation and presents some preliminaries on SHDS. Section \ref{secresults} presents the main technical assumptions and results. Section \ref{sec_examples} presents two examples, and finally Section \ref{section:conclusions} ends with the conclusions.
\section{PRELIMINARIES}
\label{sec:preli}
We denote the set of (non negative) real numbers by $(\mathbb{R}_{\geq 0})$ $\mathbb{R}$. Given a closed set $\mathcal{A} \subset \mathbb{R}^n$ and a vector $z \in \mathbb{R}^n$, we define $|z|_\mathcal{A} := \inf_{y \in \mathcal{A}}|z-y|$, and we use $|\cdot|$ to denote the standard Euclidean norm. We use $r\mathbb{B}^\circ$ to denote the open ball (in the Euclidean norm) of appropriate dimension centered around the origin and with radius $r>0$. For ease of notation, given two vectors $u,v \in \mathbb{R}^{n}$, we write $(u,v)$ for $(u^{T},v^{T})^{T}$. A function $\alpha:\mathbb{R}_{\geq0}\to\mathbb{R}_{\geq0}$ is said to be of class $\mathcal{K}$ if it is zero at zero, continuous, and strictly increasing. 
We use $\mathbf{B}(\mathbb{R}^m)$ to denote the Borel $\sigma$-field. A set $K\subset\mathbb{R}^m$ is said to be measurable if $K\in\mathbf{B}(\mathbb{R}^m)$. Given a measurable space $(\Omega, \mathcal{F})$, a set-valued map $F: \Omega \rightrightarrows \mathbb{R}^n$ is said to be \textit{$\mathcal{F}$-measurable}, if for each open set $\mathcal{O} \subset \mathbb{R}^n$, the set $F^{-1}(\mathcal{O}) :=
\{\omega \in \Omega : F(\omega) \cap \mathcal{O} 	= \varnothing\}\in \mathcal{F}$.  

\subsection{Basic Notions on Stochastic Hybrid Dynamical Systems}
\label{prelimSHDS}
In this paper, we consider SHDS that follow the formalism introduced in \cite{Teel_ANRC}. These systems are modeled as: 
\begin{subequations}\label{SHDS1}
	\begin{align}
	&y\in C,~~~~~~~~~\dot{y}= F(y),\label{SHDS_flows0}\\
	&y\in D,~~~~~~y^+= G(y,v^+),~~~v\sim \mu(\cdot),~~\label{SHDS_jumps0}
	\end{align}
\end{subequations}
where $y\in\mathbb{R}^n$ is the main state, $F:\mathbb{R}^n\to\mathbb{R}^n$ is called the flow map, $G:\mathbb{R}^n\times\mathbb{R}^m\to\mathbb{R}^n$ is called the jump map, $C$ is the flow set, $D$ is the jump set, and $v^+$ is a place holder for a sequence $\{\bf{v}_{\ell}\}_{\ell=1}^{\infty}$ of independent, identically distributed ({\em i.i.d.}) random variables ${\mathbf{v}}_k:\Omega\to\mathbb{R}^m$, $k\in\mathbb{N}$, defined on a probability space $(\Omega,\mathcal{F},\mathbb{P})$. Thus, ${\bf v_k}^{-1}(F):=\{\omega\in\Omega: \mathbf{v}_k(\omega)\in F\}\in\mathcal{F}$ for all $F\in\mathbf{B}(\mathbb{R})^m$, and $\mu:\mathbf{B}(\mathbb{R}^m)\to[0,1]$ is defined as $\mu(F):=\mathbb{P}\{\omega\in\Omega:{\bf v}_k(\omega)\in F\}$.  For the purpose of completeness, we review the concept of solution to \eqref{SHDS1}, introduced in \cite{Teel_ANRC}. Random solutions to SHDS \eqref{SHDS1} are functions of $\omega\in\Omega$, denoted ${\bf y}(\omega)$. To formally define these mappings, for $\ell\in\mathbb{Z}_{\geq1}$, let $\mathcal{F}_\ell$ denote the collection of sets $\{\omega\in\Omega:({\bf v}_1(\omega),{\bf v}_2(\omega),\ldots,{\bf {v}}_\ell(\omega))\in F\}$, $F\in\mathbf{B}(\mathbb{R}^m)^\ell)$, which are the sub-$\sigma$-fields of $\mathcal{F}$ that form the minimal filtration of ${\bf v}=\{{\bf v}_\ell \}_{\ell=1}^{\infty}$, which is the smallest $\sigma$-algebra on $(\Omega,\mathcal{F})$ that contains the pre-images of $\mathbf{B}(\mathbb{R}^m)$-measurable subsets on $\mathbb{R}^m$ for times up to $\ell$. A stochastic hybrid arc is a mapping ${\bf y}$ from $\Omega$ to the set of hybrid arcs \cite[Ch. 2]{bookHDS}, such that the set-valued mapping from $\Omega$ to $\mathbb{R}^{n+2}$, given by $\omega\mapsto \text{graph}({\bf y}(\omega)):=\big\{(t,j,z):\tilde{y}={\bf y}(\omega), (t,j)\in\text{dom}(\tilde{y}),z=\tilde{y}(t,j)\big\}$, is $\mathcal{F}$-measurable with closed-values. Let $\text{graph}({\bf y}(\omega))_{\leq \ell}:=\text{graph}({\bf y} (\omega))\cap (\mathbb{R}_{\geq0}\times\{0,1,\ldots,\ell\}\times\mathbb{R}^n).$ An $\{\mathcal{F}_\ell\}_{\ell=0}^{\infty}$ adapted stochastic hybrid arc is a stochastic hybrid arc ${\bf y}$ such that the mapping $\omega\mapsto \text{graph}({\bf y}(\omega))_{\leq \ell}$ is $\mathcal{F}_\ell$-measurable for each $\ell \in\mathbb{N}$. An adapted stochastic hybrid arc ${\bf y}(\omega)$, or simply $\mathbf{y}_\omega$,  is a solution to the SHDS \eqref{SHDS1} starting from $y_0$, denoted ${\bf y}_\omega\in \mathcal{S}_r(y_0)$, if for each $\omega\in\Omega$: (1) $\mathbf{y}_\omega(0,0)=y_0$; (2) if $(t_1,j),(t_2,j)\in\text{dom}(\mathbf{y}_{\omega})$ with $t_1<t_2$, then for almost all $t\in[t_1,t_2]$, $\mathbf{y}_{\omega}(t,j)\in C$ and $\dot{\mathbf{y}}_\omega(t,j)= F(\mathbf{y}_\omega(t,j))$; (3) if $(t,j),(t,j+1)\in\text{dom}(\mathbf{y}_\omega)$, then $\mathbf{y}_\omega(t,j)\in D$ and $\mathbf{y}_\omega(t,j+1)= G(\mathbf{y}_\omega(t,j),\mathbf{v}_{j+1}(\omega))$. We use $\mathcal{V}:=\bigcup_{\omega\in \Omega,i\in\mathbb{Z}_{\geq0}}\mathbf{v}_{i+1}$ to denote all the possible values that the random variable $v$ can take. A random solution ${\bf y}_\omega$ is said to be: a) almost surely {\em non-trivial} if its hybrid time domain contains at least two points almost surely; and b) almost surely  {\em  complete} if for almost every sample path $\omega\in \Omega$ the hybrid arc ${\bf y_\omega}$ has an unbounded time domain. If \eqref{SHDS_jumps0} does not depend on $v$, the model recovers a standard deterministic HDS \cite{bookHDS}. 
\section{MODEL AND MAIN RESULTS}
\label{secresults}
To study averaging theory in SHDS, we focus on a sub-class of systems \eqref{SHDS1}, which have continuous-time dynamics given by
\begin{subequations}\label{originalSHDS}
\begin{align}\label{flows1a}
(x,r,\tau)\in \tilde{C}:=\mathbb{R}^n\times C\times \mathbb{R}_{\geq0},~~~\left\{\begin{array}{l}
\dot{x}=f(x,r,\tau,\varepsilon)\\
\dot{r}= w(r)\\
 \dot{\tau}= \dfrac{1}{\varepsilon},
\end{array}\right.
\end{align}
and discrete-time dynamics given by
\begin{align}\label{jumps1a}
(x,r,\tau)\in\tilde{D}:=\mathbb{R}^n\times D\times\mathbb{R}_{\geq0},~~\left\{\begin{array}{l}
x^+=g(x,r,v^+)\\
r^+=h(r,v^+)\\
\tau^+=\tau
\end{array}\right.
\end{align}
\end{subequations}
where $x\in\mathbb{R}^n$ is the main state, $r\in\mathbb{R}^p$ is an auxiliary state that can be used to model timers, logic states, counters, etc. The state $\tau$ models fast time-variations in the function $f$, and its continuous-time dynamics are parameterized by a small constant $\varepsilon>0$. As in \eqref{SHDS1}, $v$ in \eqref{jumps1a} models (the same) sequence of i.i.d random inputs acting on both $g$ and $h$.

\vspace{0.1cm}
We make the following regularity assumption on the dynamics \eqref{originalSHDS}.
\vspace{0.1cm}
\begin{assumption}\label{regularityassumption1}
There exists $\varepsilon^*>0$ such that:
\begin{enumerate}
\item The function $f$ is $\mathcal{C}^1$, and the functions $w,g,h$ are continuous. The sets $C$ and $D$ are compact. 
\item The function $f$ satisfies $0=f(0,r,\tau,\varepsilon)$ for all $(r,\tau,\varepsilon)\in (C\cup D)\times\mathbb{R}_{\geq0}\times[0,\varepsilon^*)$. 
\item The function $g$ satisfies $0=g(0,r,v)$ for all $r\in C\cup D$ and all $v\in\mathcal{V}$.
\item There exists $L_x\geq0$, such that 
\begin{equation}\label{lipschitzx}
|f(x_1,r,\tau,\varepsilon)-f(x_2,r,\tau,\varepsilon)|\leq L_x|x_1-x_2|,
\end{equation}
for all $r\in C\cup D$, $x_1,x_2\in\mathbb{R}^n$, $\tau\in\mathbb{R}_{\geq0}$, and $\varepsilon\in[0,\varepsilon^*)$.
\item There exists $L_{\varepsilon^*}\geq0$, such that
\begin{equation}\label{lipschitzep}
|f(x,r,\tau,\varepsilon_1)-f(x,r,\tau,\varepsilon_2)|\leq L_{\varepsilon^*}|x||\varepsilon_1-\varepsilon_2|
\end{equation}
for all $r\in C\cup D$, $x\in\mathbb{R}^n$, $\tau\in\mathbb{R}_{\geq0}$, and $\varepsilon_1,\varepsilon_2\in(0,\varepsilon^*)$.
\item There exists $L_g>0$ such that
\begin{equation}\label{Lipscthizconditionjumps}
|g(x_1,r,v)-g(x_2,r,v)|\leq L_g|x_1-x_2|,
\end{equation}
for all $r\in C\cup D$, all $x_1,x_2\in\mathbb{R}^n$, and all $v\in \mathcal{V}$.
\item The map $$G(x,r,v):=g(x,r,v)\times h(r,v)$$ satisfies $G(\mathbb{R}^n\times D\times\mathcal{V})\subset \mathbb{R}^n\times (C\cup D)$, and $|h(r,v)|\leq H$, for some $H>0$ and all $(r,v)\in D\times\mathcal{V}$. 
\end{enumerate}
\end{assumption}

\vspace{0.1cm}
\begin{remark}
To contextualize Assumption \ref{regularityassumption1}, we note that the conditions imposed by items 2), 4), and 5), are the same as those typically studied for smooth ODEs \cite[pp. 170]{Sastry:89} . In our case, since we are interested in global results, we ask that these conditions hold globally. On the other hand, conditions 1), 3), and 6) are reasonable extensions to the discrete-time case. Note that the uniform Lipschitz assumption over $r$ is not restrictive because $C\cup D$ is compact. Finally, condition 7) simply asks that jumps always map back to $C\cup D$. \hfill \QEDB
\end{remark}

\vspace{0.1cm}
\begin{remark}
Assumption \ref{regularityassumption1} rules out finite escape times in \eqref{flows1a}, which is needed in order to use Lyapunov-Foster functions to certify stability properties of unbounded closed sets in SHDS. These types of sets will necessarily emerge in \eqref{originalSHDS} because $\tau$ will grow unbounded. \QEDB
\end{remark}

\vspace{0.1cm}
Next, we introduce the average map of $x$ in \eqref{flows1a}.

\vspace{0.1cm}
\begin{definition}\label{definitionaveragemap}
The function $f(\cdot,\cdot,\cdot,0)$ is said to have an \emph{average map} $f_{\text{ave}}(\cdot)$ if there exists a continuous, bounded, and non-increasing function $\gamma:\mathbb{R}_{\geq0}\to\mathbb{R}_{\geq0}$ satisfying $\lim_{\tau\to\infty}\gamma(\tau)=0$, such that
\begin{equation}\label{average0assumption}
\Big|\frac{1}{T}\int_{\tau}^{\tau+T} \big(f(x,r,s,0)-f_{\text{ave}}(x,r)\big)\text{d}s\Big|\leq \gamma(T)|x|, 
\end{equation}
for all $x\in\mathbb{R}^n$, all $r\in C\cup D$, and all $T,\tau\in\mathbb{R}_{\geq0}$. The function $\gamma(\cdot)$ is called the \emph{convergence function}. \QEDB
\end{definition}

\vspace{0.1cm}
To simplify notation, we now introduce the function
\begin{equation}\label{distanceaverage00}
d(x,r,\tau):=f(x,r,\tau,0)-f_{\text{ave}}(x,r),
\end{equation}
which will also play a prominent role in our analysis.  As in Assumption \ref{regularityassumption1}, we will assume that $d$ and $f_{\text{ave}}$ have suitable smoothness properties on $\mathbb{R}^n$, uniformly on $C\cup D$.

\vspace{0.1cm}
\begin{assumption}\label{smootaverageassumption}
The function $f(\cdot,\cdot,\cdot,0)$ admits a $\mathcal{C}^1$ average map $f_{\text{ave}}$ with convergence function $\gamma$, and there exists $L_{\text{ave}}\geq0$ such that:
\begin{equation}
|f_{\text{ave}}(x_1,r)-f_{\text{ave}}(x_2,r)|\leq L_{\text{ave}}|x_1-x_2|,
\end{equation}
for all $x_1,x_2\in\mathbb{R}^n$, and all $r\in C\cup D$. Moreover,
\begin{equation}\label{average200}
\Big|\frac{1}{T}\int_{\tau}^{\tau+T} \frac{\partial d(x,r,\tau)}{\partial(x,r)} \text{d}\tau\Big|\leq \gamma(T).
\end{equation}
for all $(x,r,\tau)\in\mathbb{R}^n\times (C\cup D)\times\mathbb{R}_{\geq0}$ and $T\geq0$. \QEDB
\end{assumption}

\vspace{0.1cm}
\begin{remark}
Inequality \eqref{average0assumption}, used in Definition \ref{definitionaveragemap}, and the conditions of Assumption \ref{smootaverageassumption}, are equivalent to the assumptions made in the analysis of smooth ODEs based on averaging theory  \cite[Lemma 4.2.2]{Sastry:89}. Therefore, our assumptions are not more restrictive than those considered in the literature for ODEs, again, subject to the fact that we aim for global results as opposed to local or semi-global. In fact, if the state $x$ in \eqref{originalSHDS} is restricted to evolve in a compact set, then most of the global Lipschitz conditions can be relaxed.
 \QEDB
\end{remark}

\vspace{0.1cm}
\begin{remark}
An important class of systems that satisfy Assumptions \ref{regularityassumption1}-\ref{smootaverageassumption}  are linear time-varying systems, affine in $r$, with matrices $A(t)$ generating a well-defined average matrix $\bar{A}$. These systems usually emerge in adaptive control and parameter estimation via least squares, see \cite[Sec. 4.1]{Sastry:89}. Similar algorithms that make use of stochastic resets were also recently studied in \cite{OchoaEstimation2023}. \QEDB
\end{remark}

\vspace{0.1cm}
Next, we introduce the Average SHDS of system \eqref{originalSHDS}.

\vspace{0.1cm}
\begin{definition}\label{averagesystem}
The Average SHDS of system \eqref{originalSHDS} has state $z=(\hat{x},\hat{r})\in\mathbb{R}^n\times (C\cup D)$, and dynamics
\begin{subequations}\label{averagehybridynamics}
\begin{align}
&z\in\mathbb{R}^n\times C,~~~~\dot{z}= F_{\text{ave}}(z):=\left(\begin{array}{c}
f_{\text{ave}}(\hat{x},\hat{r})\\
w(\hat{r})
\end{array}\right)\label{flowsaveragedefinition}\\
&z\in\mathbb{R}^n\times D,~~z^+=G_{\text{ave}}(z,v^+),\label{jumpsaveragedefinition}
\end{align}
\end{subequations}
where the average jump map is defined as $G_{\text{ave}}(z,v):=g(z,v)\times h(\hat{r},v)$, for all $(z,v)\in(\mathbb{R}^n\times D)\times \mathcal{V}$. \QEDB
\end{definition}

\vspace{0.1cm}
We can now state the following Lemma, which is a straightforward consequence of Assumptions \ref{regularityassumption1}-\ref{smootaverageassumption}, and \cite{Teel_ANRC}.

\vspace{0.1cm}
\begin{lemma}
The SHDS \eqref{originalSHDS} and the Average SHDS \eqref{averagehybridynamics} satisfy the Stochastic Hybrid Basic Conditions of \cite[Assumptions 1 $\&$ 2]{Teel_ANRC}. \QEDB
\end{lemma}

\vspace{0.1cm}

Our goal is to leverage the stability properties of the Average SHDS \eqref{averagehybridynamics}  to study the stability properties of the original SHDS \eqref{originalSHDS}. To achieve this task, we make use of the following assumption, which is a stochastic hybrid extension of the assumptions made in \cite[Thm. 4.2.5]{Sastry:89} for ODEs.

\vspace{0.1cm}
\begin{assumption}[Stability via Lyapunov-Foster Function]\label{Lyapunovassumptionaverage}
Let $\mathcal{A}:=\{0\}\times (C\cup D)$. There exists a smooth function $V:\mathbb{R}^n\times (C\cup D) \to\mathbb{R}_{\geq0}$, and constants $c_i>0$, such that:
\begin{align}\label{inequalitiesLyapunov}
&c_1|z|_{\mathcal{A}}^2\leq V(z)\leq c_2|z|_{\mathcal{A}}^2,~~~|\nabla V(z)|\leq c_3|z|_{\mathcal{A}},
\end{align}
for all $z\in\mathbb{R}^n\times (C\cup D)$. Similarly, 
\begin{equation}
\langle \nabla V(z),F_{\text{ave}}(z)\rangle \leq -c_4 V(z),~\forall~z\in\mathbb{R}^n\times C,
\end{equation}
and the expected value of $V$ during jumps of \eqref{jumpsaveragedefinition} satisfies:
\begin{align}\label{Lyapunovjumpsstochastic}
\hspace{-0.6cm}\int_{\mathbb{R}^m}V(G_{\text{ave}}(z,v))\mu(dv)\leq c_5V(z),~\forall~z\in\mathbb{R}^n\times D,
\end{align}
where $\lambda:=\frac{c_2}{c_1}c_5<\frac{1}{2}$. \QEDB
\end{assumption}
\vspace{0.2cm}
\begin{remark}
By directly invoking \cite[Thm. 4.2]{Teel_ANRC}, Assumption \ref{Lyapunovassumptionaverage} guarantees that the compact set $\mathcal{A}$ is Uniformly Globally Asymptotically Stable in Probability (UGASp) for system \eqref{averagehybridynamics} \cite[Sec. 2.3]{Teel_ANRC}. In fact, by  \cite[Thm. 11]{teel2015stochastic}, the set $\mathcal{A}$ is also Uniformly Globally Exponentially Stable in the Mean (UGES-M) \cite[pp. 3129]{teel2015stochastic}, i.e., there exist  $k_1, k_2 >0$  such that every random solution $\mathbf{z}_{\omega}$ of \eqref{averagehybridynamics} satisfies: 
\begin{equation}
\mathbb{E}\left[|\mathbf{z}_{\omega}(\textbf{T},\textbf{J})|_{\mathcal{A}}e^{k_2 (\textbf{T}+\textbf{J})}\right]\leq k_1 |z(0,0)|_{\mathcal{A}},
\end{equation}
for all hybrid stopping times $(\textbf{T},\textbf{J})$ of $\mathbf{z}_{\omega}$ \cite[pp. 3126]{teel2015stochastic}, and for all $z(0,0)\in\mathbb{R}^n$. \QEDB
\end{remark}

\vspace{0.1cm}
In some cases, characterizing a stable set for a SHDS might be challenging or even impossible if such set does not exist. In this case, the property of \emph{recurrence} \cite[Sec. 2.4]{teel2015stochastic}, can still be studied via Foster functions. 
\vspace{0.05cm}
\begin{definition}\label{definition3}
An open (relative to the flow and jump sets), bounded set $\mathcal{O}\subset\mathbb{R}^{n}$ is \emph{uniformly globally recurrent} (UGR) for a SHDS \eqref{SHDS1} if there are no finite escape times for \eqref{SHDS_flows0}, and for each $\rho>0$ and $R>0$ there exists $\tau\geq0$ such that every maximal random solution $\bf{y}_{\omega}$ starting in $R\mathbb{B}$ satisfies

\vspace{-0.3cm}
\begin{small}
\begin{align*}
&\mathbb{P}\Big(\big(\text{graph}({\bf y}_{\bf{\omega}})\subset (\Gamma_{<\tau}\times\mathbb{R}^n)\big)\lor \big(\text{graph}({\bf y}_{\bf\omega})\cap (\Gamma_{\leq\tau}\times\mathcal{O})\big)\Big)\\
&~~~~~~~~~~~~~~~~~~~~~~~~~~~~~~~~~~~~~~~~~~~~~~~~~~~~~~~~~~\geq 1-\rho,
\end{align*}
\end{small}

\vspace{-0.2cm}\noindent 
where $\text{graph}({\bf y_\omega}):=\big\{(t,j,s):s={\bf y}_\omega, (t,j)\in\text{dom}(y)\big\}$, and $\Gamma_{<\tau}:=\{(s,t)\in\mathbb{R}^2:s+t<\tau\}$.  \QEDB
\end{definition}

\vspace{0.1cm}
Loosely speaking, Definition \ref{definition3} says that from every initial condition, solutions to \eqref{SHDS1} either stop or hit the set $\mathcal{O}$, with a hitting time that is uniform over compact sets of initial conditions, and the solutions do not have finite escape times.
\begin{remark}
The property of recurrence is related to notions of  ``metastability'', which emerge in walking robots \cite{Metastability,BylWalking} and UAVs performing stochastic exploration \cite{ACC_recurrence}. \QEDB
\end{remark}

\vspace{0.1cm}
Next, we present the main result of this paper. For the original SHDS \eqref{originalSHDS}, we establish stability, or recurrence, of sets equal or related to $\mathcal{A}\times\mathbb{R}_{\geq0}$, based on the stability properties of $\mathcal{A}$ for the Average SHDS \eqref{averagehybridynamics}.
\vspace{0.15cm}

\begin{theorem}\label{maintheorem}
Suppose that Assumptions \ref{regularityassumption1}-\ref{Lyapunovassumptionaverage} hold. Then, there exists $\varepsilon^\star>0$ and a class $\mathcal{K}$ function $\varphi$ such that for all $\varepsilon\in(0,\varepsilon^\star)$:
\begin{enumerate}
\item The set $(\mathcal{A}+\varphi(\varepsilon)\mathbb{B}^{\circ})\times \mathbb{R}_{\geq0}$ is UGR for the SHDS \eqref{originalSHDS}.
\item If, additionally, $\frac{\partial f(x,r,\tau,\varepsilon)}{\partial r}\dot{r}=0$ for all $(x,r,\tau,\varepsilon)\in\mathbb{R}^n\times C\times \mathbb{R}_{\geq0}\times[0,\varepsilon^*]$, then the set $\mathcal{A}\times \mathbb{R}_{\geq0}$ is UGES-M for the SHDS \eqref{originalSHDS}.
\end{enumerate}
\QEDB
\end{theorem}
The proof of Theorem \ref{maintheorem} is long, and it will be presented on a different occasion.

\section{EXAMPLES}
\label{sec_examples}

\vspace{-0.1cm}
We present two simple examples to illustrate Theorem \ref{maintheorem}. In both cases, a continuous-time plant or controller, with main state $x\in\mathbb{R}$, and high-frequency terms acting on the flow map, is jammed every $T$ seconds according to $x^+=(0.75+v)x$ , where $v\in\{-0.75,0.75\}$ is a (scaled) Bernoulli variable that satisfies $v=0.75$ with probability $p$. This behavior can be modeled using a periodic resetting timer $r$, leading to the SHDS with continuous-time dynamics
\begin{subequations}\label{examples1}
\begin{align}\label{example1F}
(x,r,\tau)\in \mathbb{R}\times [0,T]\times \mathbb{R}_{\geq0},~\left\{\begin{array}{l}
\dot{x}=f(x,\tau)\\
\dot{r}= 1\\
 \dot{\tau}= \dfrac{1}{\varepsilon},
\end{array}\right.
\end{align}
and discrete-time dynamics
\begin{align}\label{example1J}
(x,r,\tau)\in\mathbb{R}\times \{T\}\times\mathbb{R}_{\geq0},~\left\{\begin{array}{l}
x^+=(0.75+v)x\\
r^+=0\\
\tau^+=\tau
\end{array}\right.
\end{align}
\end{subequations}
Since $f$ depends on the fast-varying signal $\tau$, we are interested in characterizing the resilience of the system with respect to both the high-frequency variation of $f$, and the random jamming. To do this, we will study the behavior of the trajectories $(x,r)$ with respect to the compact set $\mathcal{A}=\{0\}\times[0,T]$.
\vspace{0.1cm}
\subsubsection{Linear Actuator with State-Dependent Disturbance} We first consider the case when $x$ models the state of an actuator with dynamics
$$\dot{x}=-x(1+\sin(\tau))+u,$$ where $u$ is a constant input and $\sin(\tau)$ is the high-frequency state-dependent disturbance, see also \cite[Ch. 10]{KhalilBook} To study stability of the origin, we set $u=0$. Clearly, the function $f(x,\tau)=-x(1+\sin(\tau))$ satisfies Assumption \ref{regularityassumption1}. The average dynamics of \eqref{example1F}-\eqref{example1J} are given by a SHDS with continuous-time dynamics:
\begin{equation}\label{exampleaverage11}
\hspace{-1.1cm}(\hat{x},\hat{r})\in \mathbb{R}\times [0,T],~~~\left\{\begin{array}{l}
\dot{\hat{x}}=-\hat{x}\\
\dot{\hat{r}}= 1
\end{array}\right.
\end{equation}
and discrete-time dynamics
\begin{equation}\label{exampleaverage12}
~~~~~~~~(\hat{x},\hat{r})\in\mathbb{R}\times \{T\},~~~~\left\{\begin{array}{l}
\hat{x}^+=(0.75+v)\hat{x}\\
\hat{r}^+=0,
\end{array}\right.
\end{equation}

\vspace{-0.1cm}\noindent
which also satisfies Assumption \ref{smootaverageassumption}. Assumption \ref{Lyapunovassumptionaverage} can be verified using the Lyapunov-Foster function $V(\hat{x})=\hat{x}^2$. In this case, \eqref{inequalitiesLyapunov} holds with $c_1=c_2=1$, $c_3=c_4=2$. Moreover, $$\mathbb{E}[V(\hat{x}^+)]=p(1.5\hat{x})^2+(1-p)(0\hat{x})^2=\frac{9}{4} pV(\hat{x}),$$ which implies that \eqref{Lyapunovjumpsstochastic} is satisfied with $c_5>0$ and $\lambda\in(0,\frac{1}{2})$ if $p<\frac{2}{9}\approx0.222$. Since in this case $\frac{\partial f(x,r,\tau,\varepsilon)}{\partial r}\dot{r}=0$, we can directly apply Theorem \ref{maintheorem} to conclude that, provided the probability of jamming satisfies $p\in(0,0.22)$, there exists $\varepsilon^*>0$, such that for all $\varepsilon\in(0,\varepsilon^*)$, the set $\mathcal{A}\times\mathbb{R}_{\geq0}$ is UGES-M for the SHDS \eqref{example1F}-\eqref{example1J}.

\vspace{0.1cm}
\subsubsection{Global Optimization Under Random Jamming} 
\begin{figure}[t!]
\centering 
\includegraphics[width=0.99\columnwidth]{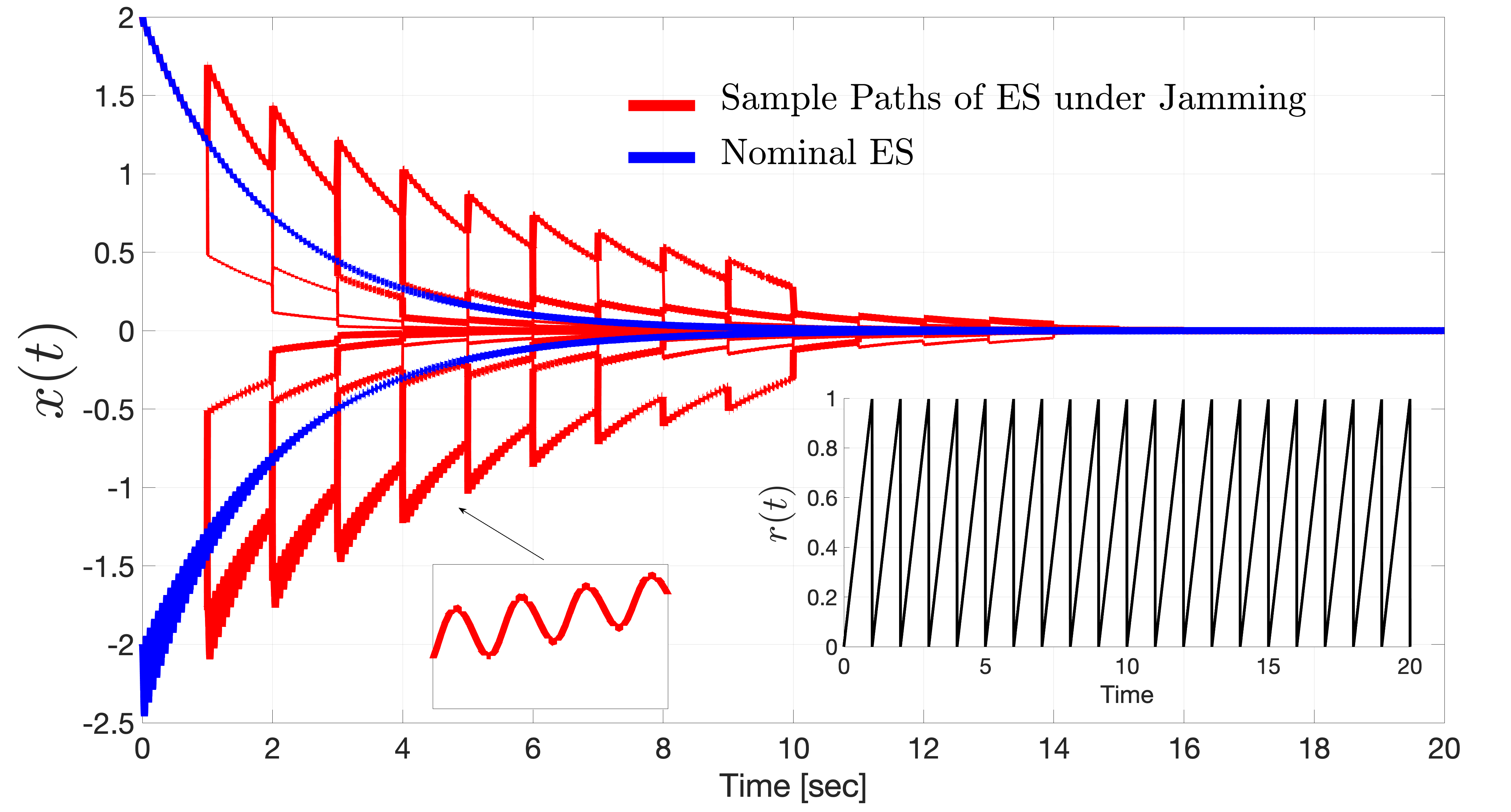}
\caption{Multiple sample paths generated by system \eqref{examples1}.}
\label{fig:toyexample2}
\vspace{-0.3cm}
\end{figure}

Consider now the case when the function $f(x,\tau)$ in \eqref{example1F} is
\begin{equation}
f(x,\tau)=-\frac{1}{a}\phi(x+a\sin(\tau))\sin(\tau),
\end{equation}
where $a>0$ is a tunable parameter.
This function describes the dynamics of a standard extremum seeking controller \cite{DerivativesESC,KrsticBookESC,PoTe17Auto} aiming to minimize a cost function $\phi:\mathbb{R}\to\mathbb{R}$, using a sinusoidal dither. To illustrate the resilience properties of these dynamics with respect to random jamming, we consider the cost $\phi(y)=y^2$,  leading to 
$$f(x,\tau)=-\frac{x^2}{a}\sin(\tau)-2x\sin(\tau)^2-a\sin(\tau)^3.$$
This function does not satisfy Assumption \ref{regularityassumption1}. However, we can design the dither amplitude so that Assumption \ref{regularityassumption1} holds outside a small neighborhood of the origin. Indeed, using\footnote{Naturally, such a choice requires knowledge of $\phi$, and wouldn't be feasible in model-free applications. We use it here only for the purpose of illustration.} $a(x)=\max\{\delta,|x|\}$, where $\delta>0$ is a small constant, leads now to the following dynamics in $\mathbb{R}\backslash\delta\mathbb{B}$ for the state $x$:
\begin{equation}\label{ESglobal}
\dot{x}=f(x,\tau)=-x\sin(\tau)-2x\sin(\tau)^2-|x|\sin(\tau)^3.
\end{equation}
which satisfies Assumption \ref{regularityassumption1} for all $x\in\mathbb{R}\backslash\delta\mathbb{B}$. Moreover, since the average of $\sin(\tau)^2$ is $1/2$, and the average of $\sin(\tau)^3=0$, the resulting SHDS \eqref{example1F}-\eqref{example1J} has the exact same average system \eqref{exampleaverage11}-\eqref{exampleaverage12}, for which the same Lyapunov-Foster function can be used to study the stability properties of the set $\mathcal{A}$. It follows that all the main bounds of our proofs hold in $\mathbb{R}\backslash\delta\mathbb{B}$. Since in $\delta\mathbb{B}$ the dynamics of $x$ are locally Lipschitz, the SHDS still satisfies the Stochastic Hybrid Basic Conditions, and we can conclude global recurrence of the set $(\mathcal{A}+\delta\mathbb{B}^{\circ})\times\mathbb{R}_{\geq0}$ for system \eqref{examples1}. Figure \ref{fig:toyexample2} presents 100 sample paths of \eqref{examples1}  simulated from initial conditions satisfying $x(0,0)\in\{-2,2\}$, using the function $f(x,\tau)$ defined in \eqref{ESglobal}. We also show, in the color blue, a trajectory generated by the nominal dynamics \eqref{ESglobal} without jamming. The inset shows the evolution in time of the resetting clock $r$, indicating the times when the state $x$ was randomly jammed. As observed, on average, the ES dynamics under jamming still converge to a neighborhood of the origin.

\section{CONCLUSIONS}
\label{section:conclusions}
An averaging result was introduced for a class of stochastic hybrid dynamical systems with time-varying flow maps having a well-defined average mapping. By using Lyapunov-Foster functions (and a suitable state transformation), global stability or recurrence can be established for the original dynamics based on the stability properties of the average stochastic hybrid dynamics. Two simple examples were presented to illustrate the main result. The assumptions considered in the paper parallel those considered in the literature of ODEs whenever global results (as opposed to semi-global) are sought-after. However, some of our assumptions can be further relaxed to provide more flexibility in the analysis and design of stochastic hybrid algorithms and controllers via averaging theory. Of particular interest are SHDS with set-valued maps,  as well as applications to stochastic hybrid extremum seeking control and hybrid vibrational control. 

\bibliographystyle{IEEEtran}
\bibliography{dissertationA.bib}

\end{document}